\documentclass{dmtcs}

\usepackage[latin1]{inputenc}
\usepackage{subfigure}

\usepackage{amssymb}
\usepackage{amsmath}
\usepackage{bbm}
\newcommand{\be}{\[}
\newcommand{\ee}{\]}
\newcommand{\ds}{\displaystyle}
\newcommand{\id}[1]{\mathbbm{1}_{#1}}

\newcommand{\bs}[1]{\boldsymbol{#1}}
\newcommand{\bra}[1]{\langle \, #1 \,|}
\newcommand{\ket}[1]{|\, #1 \,\rangle}

\newtheorem{thm}{Theorem}
\newtheorem{prop}[thm]{Proposition}
\newtheorem{cor}[thm]{Corollary}
\newtheorem{lem}[thm]{Lemma}
\newtheorem{defn}{Definition}

\newtheorem{conj}{Conjecture}

\author{Arvind Ayyer\addressmark{1}\and
  Volker Strehl\addressmark{2}} 
\title[Asymmetric annihilation process]{The spectrum of an asymmetric annihilation process} 
\address{\addressmark{1}
Institut  de Physique Th\'eorique, C. E. A.  Saclay,
 91191 Gif-sur-Yvette Cedex, France \\
\addressmark{2}
Department of Computer Science, Universit\"at Erlangen-N\"urnberg,
Haberstrasse 2, D-91058 Erlangen, Germany } 
\begin{document}
\maketitle
\begin{abstract}
\paragraph{Abstract.}
In recent work on nonequilibrium statistical physics, a certain
Markovian exclusion model called an asymmetric annihilation process
was studied by Ayyer and Mallick.  In it they gave a precise
conjecture for the eigenvalues (along with the multiplicities) of the
transition matrix. They further conjectured that to each eigenvalue, there
corresponds only one eigenvector.
We prove the first of these conjectures by generalizing the original
Markov matrix by introducing extra parameters, explicitly calculating
its eigenvalues, and showing that the new matrix reduces to the
original one by a suitable specialization.  
In addition, we outline a
derivation of the partition function in the generalized model, which
also reduces to the one obtained by Ayyer and Mallick in the original
model.

\paragraph{R\'esum\'e.}
Dans un travail r\'ecent sur la physique statistique hors \'equilibre, un
certain mod\`ele d'exclusion Markovien appel\'e ``processus d'annihilation
asym\'etriques'' a \'et\'e \'etudi\'e par Ayyer et Mallick. Dans ce document,
ils ont donn\'e une conjecture pr\'ecise pour les valeurs propres (avec
les multiplicit\'es) de la matrice stochastique. Ils ont en outre suppos\'e
que, pour chaque valeur propre, correspond un seul vecteur propre.
Nous prouvons la premi\`ere de ces conjectures en g\'en\'eralisant la
matrice originale de Markov par l'introduction de param\`etres
suppl\'ementaires, calculant explicitement ses valeurs propres, et en
montrant que la nouvelle matrice se r\'eduit à l'originale par une
sp\'ecialisation appropri\'ee.  En outre, nous pr\'esentons un calcul de la
fonction de partition dans le mod\`ele g\'en\'eralis\'e, ce qui r\'eduit
\'egalement \`a celle obtenue par Ayyer et Mallick dans le mod\`ele original.

\end{abstract}

  \section{Introduction} 
  \label{sec:intro}
In the past few years, special stochastic models motivated by
nonequilibrium statistical mechanics have motivated several
combinatorial problems. The most widely studied problem among these
has been the totally asymmetric simple exclusion process (TASEP).  The
model is defined on a one dimensional lattice of $L$ sites, each site
of which either contains a particle or not.  Particles in the interior
try to jump with rate 1 to the site to the right. The jump succeeds if
that site is empty and fails if not. On the boundary, particles enter
with rate $\alpha$ on the first site if it is empty and leave from the
last site with rate $\beta$.  This was first solved  in 1993 by developing a
new technique now called the {\em matrix product representation} \cite{DEHP}.

It was initially studied in a combinatorial setting by Shapiro and
Zeilberger in an almost forgotten paper \cite{sz} in 1982, but only
after the steady state distribution of the model was explicitly
presented in \cite{DEHP}, the problem gained widespread attention. One
of the reasons for this interest was that the common denominator of
the steady state probabilities for a system of size $L$ was $C_{L+1}$,
the $(L+1)$-th Catalan number.
One of the first articles to explain this fact combinatorially was the
one by Duchi and Schaeffer \cite{ds}, who enlarged the space of
configurations to one in bijection with bicolored Motzkin paths and
showed that the steady state distribution was uniform on this
space. The analogous construction for the partially asymmetric version of
the model (PASEP) has been done in \cite{bcepr}.

Further work has been on the relationship of the total and partially
asymmetric  exclusion processes to
different kinds of tableaux by Corteel and Williams \cite{cw1, cw2, cw3}
(permutation tableaux, staircase tableaux) and by Viennot \cite{v1}
(Catalan tableaux), to lattice paths \cite{cjpr}, and to Askey-Wilson
polynomials \cite{cw2}.

Just like the common denominator for the TASEP of size $L$ was the
Catalan number $C_L$ (which has many combinatorial interpretations),
the common denominator for the asymmetric annihilation process
considered in \cite{am} in a system of size $L$ at
$\alpha=1/2,\beta=1$ is $2^{\binom{L+1}2}$ which is the number of
domino tilings of an Aztec diamond of size $L$ as well as the number
of 2-enumerated $L \times L$ alternating sign matrices. One can
therefore hope to enlarge the configuration space as was done for the
TASEP \cite{ds} to explain this phenomena.

The remainder of this extended abstract is organized as follows: 
In Sec. \ref{sec:mod} we describe the model of the asymmetric annihilation process.
In Sec. \ref{sec:tma} we present some of the main results obtained by Ayyer and Mallick in \cite{am}.
Their work lead to a conjecture about the spectrum of this process. 
In Sec. \ref{sec:mm} we prove this conjecture by appropriately extending the model
and viewing it in a different basis obtained by a variant of the Hadamard transform.
In the concluding section we outline the derivation of the partition function for
the generalized model using the same transformation, an approach very different from
the way Ayyer and Mallick obtained the partition function in the original model.

\section{The model}
 \label{sec:mod}

Motivated by Glauber dynamics of the Ising model, Ayyer and Mallick
\cite{am} considered a 
non-equilibrium system on a finite lattice with $L$ sites labelled from
1 to $L$. States of the system are encoded by bitvectors 
$\bs{b}=b_1b_2\ldots b_L$ of length $L$, where 
$b_j \in \mathbb{B}=\{0,1\}$, so that we have a total of $2^L$ states.
These bit vectors may  be represented numerically using the
binary expansion 
$
	(\bs{b})_2 = b_L+b_{L-1}\,2^1+b_{L-2}\,2^2 + \cdots + b_1\,2 ^{L-1},
$
which introduces a total order on $\mathbb{B}^L$, so that we shall write
$\bs{b} < \bs{c}$ iff $(\bs{b})_2 < (\bs{c})_2$. All  matrices and vectors are
indexed w.r.t. this order.

The evolution rules of the system introduced in \cite{am} can now be stated as
rewrite rules for bit vectors:

\begin{itemize}
\item  In the bulk we have \emph{right shift} and \emph{annihilation} given by
\be \label{bulk1}
\begin{array}{lcl}
\text{right shift}		& 10 \to 01	& \text{with rate}~1, \cr
\text{annihilation} 	& 11 \to 00 	& \text{with rate}~\lambda,
\end{array}
\ee
and visualized (for $L=8$) in Fig. \ref{fig:shift}.
\begin{figure}[htbp]  %  figure placement: here, top, bottom, or page
   \centering
   \includegraphics[width=2in]{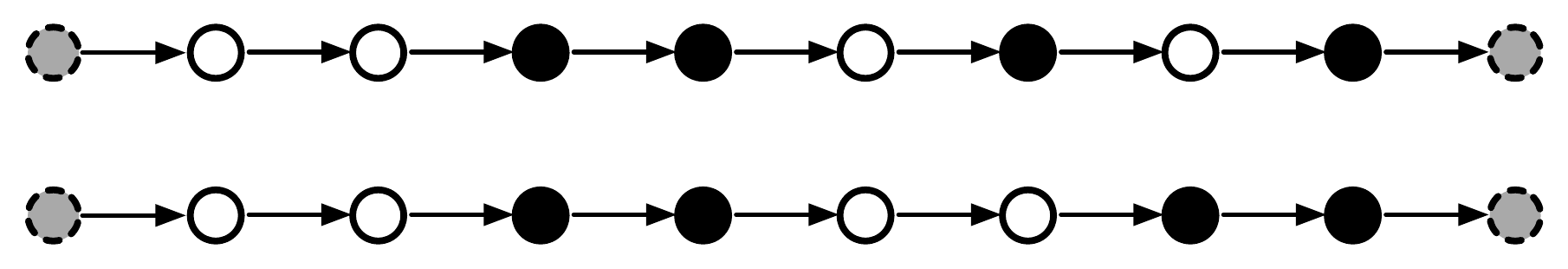} ~~\hskip1cm
   \includegraphics[width=2in]{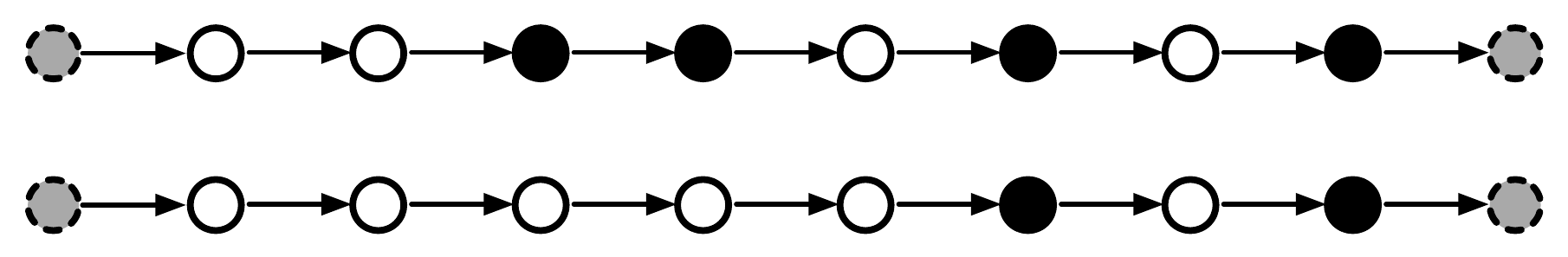}
   \caption{Right shift $00110\underline{10}1 \to 00110\underline{01}1$ and annihilation 
   		$00\underline{11}0101 \to 00\underline{00}0101$}
   \label{fig:shift}
\end{figure}

\item
On the left boundary, particles enter by \emph{left creation}
in a way consistent with the
bulk dynamics. A particle at site 1 may also be \emph{left annihilated} (due to a virtual
particle at site 0). Therefore, the first  site evolves as
\be \label{lbound}
\begin{array}{lcl}
\text{left creation}		& 0 \to 1	& \text{with rate}~\alpha, \cr
\text{left annihilation} 	& 1 \to 0 	& \text{with rate}~\alpha\lambda,
\end{array}
\ee
as illustrated Fig. \ref{fig:left}.
\begin{figure}[htbp]  %  figure placement: here, top, bottom, or page
   \centering
   \includegraphics[width=2in]{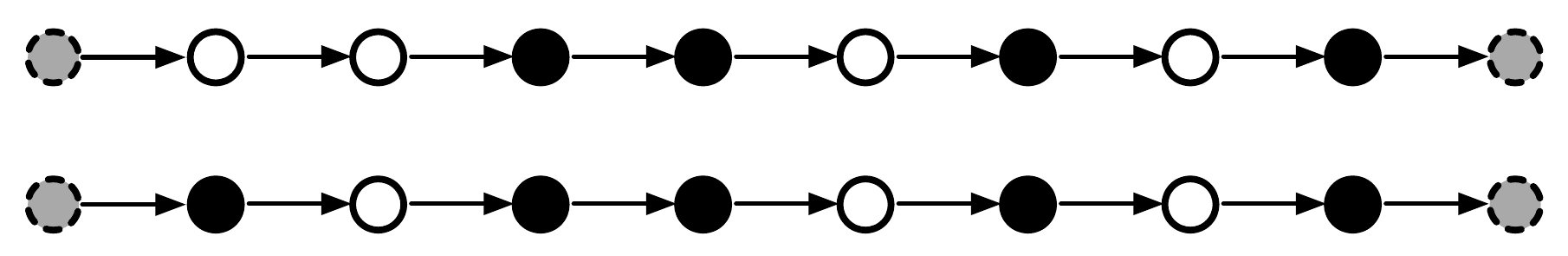} \hskip1cm
    \includegraphics[width=2in]{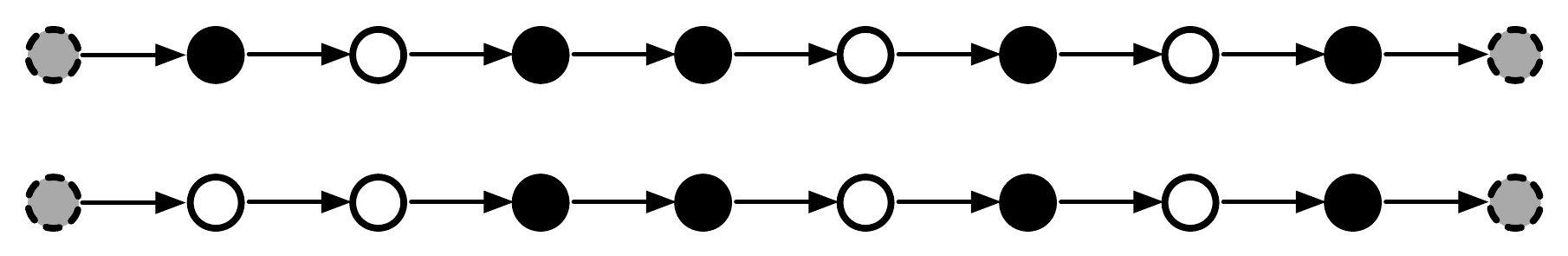}  
   \caption{Left creation $\underline{0}0110101 \to \underline{1}0110011$ and left annihilation 
   		$\underline{1}0110101 \to \underline{0}0110101$}
   \label{fig:left}
\end{figure}
\item
Particles can exit from the last site by \emph{right annihilation} 
(with a virtual particle at site $L+1$) according to 
\be \label{rbound}
\begin{array}{lcl}
\text{right annihilation} 	& 1 \to 0 	& \text{with rate}~\beta,
\end{array}
\ee
as illustrated by Fig. \ref{fig:right}.
\begin{figure}[htbp] %  figure placement: here, top, bottom, or page
   \centering
   \includegraphics[width=2in]{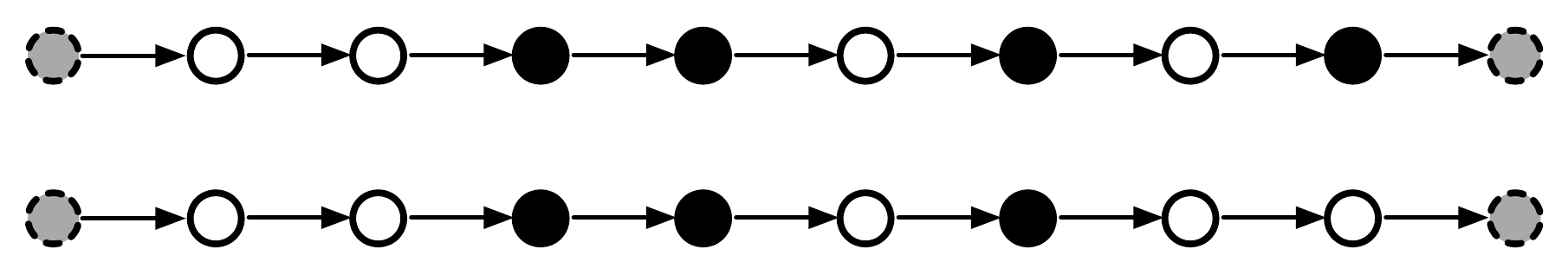}  
   \caption{Right annihilation $0011010\underline{1} \to 1011001\underline{0}$} 
   \label{fig:right}
\end{figure}
\end{itemize}
Note that all transition rules except left creation are monotonically decreasing
w.r.t. the natural order of bit vectors. Thus the transition matrix, as discussed
in the next section, is not in triangular shape.

Following \cite{am}, we will take $\lambda=1$ as that is the
only case for which they derive explicit formulae.

\section{Algebraic properties of the model} \label{sec:tma}

Is this section we present without proofs the main results as obtained by
Ayyer and Mallick in \cite{am}. First recall the general concept:

\begin{defn}
A (continuous-time) {\em transition matrix} or {\em Markov matrix} or {\em stochastic
matrix} is a square matrix of size equal to the cardinality of the
configuration space whose $(i,j)$-th entry is given by the rate of the
transition from configuration $j$ to configuration $i$, when $i$ is
not equal to $j$. The $(i,i)$-th entry is then fixed by demanding that
the entries in each column sum to zero.
\end{defn}

The Markov chain we defined in the previous section
 satisfies what \cite{am} call the ``transfer matrix Ansatz''.
The following general definition applies 
to any family of Markov processes defined by Markov matrices $\{M_L \}$ 
of increasing sizes (in most physical applications, $L$ is the size of the system).  

\begin{defn}
A  family $M_L$ of  Markov processes  satisfies the {\em Transfer Matrix
  Ansatz} if there exist matrices $T_{L,L+1}$ for all sizes $L$
such that 
\be \label{tm}
\mathrm{(TMA)}~~~~M_{L+1} T_{L,L+1} = T_{L,L+1} M_L \, .
\ee
 We also  impose that  this equality is nontrivial in the sense that
$
M_{L+1} T_{L,L+1} \neq 0. \label{cond:nontriv}
$
\end{defn}
   The  rectangular transfer matrices $T_{L,L+1}$ can be interpreted as {\em
  semi-similarity transformations} connecting Markov matrices of different
sizes.

The  last  condition is important because there is always a trivial
solution whenever we are guaranteed a unique Perron-Frobenius
eigenvector for all transition matrices $M_L$.
If  $|v_L\rangle$
 is this eigenvector of $M_L$ and 
 $\langle 1_L| = (1,1,\ldots,1)$,  
   the matrix $V_{L,L+1} =
|v_{L+1}\rangle\langle 1_L| $ satisfies (TMA)
since the Markov matrices satisfy the conditions
$\langle 1_L| M_L = 0$ and $M_{L+1}|v_{L+1} \rangle=0$.

 The above definition leads immediately to a recursive computation of
 the steady state vector, which is the zero eigenvector.  First we have
\be
 \label{constructvL}
0 = T_{L,L+1} M_L |v_L \rangle = M_{L+1} T_{L,L+1} |v_L \rangle,
\ee
which, assuming $T_{L,L+1} |v_L \rangle \neq 0$, and taking
 into account the uniqueness of
the steady state, allows us to define $|v_{L+1} \rangle$ so that 
\be \label{nontriv}
T_{L,L+1} |v_L \rangle = |v_{L+1} \rangle.
\ee
This is very analogous to the matrix product representation of
\cite{DEHP} because the steady state probability of any configuration
of length $L+1$ is expressed as a linear combination of those of
length $L$. The transfer matrix Ansatz is a stronger requirement than
the matrix product representation in the sense that not every system
which admits the representation satisfies the Ansatz. For example, the
only solution for 
(TMA) in the case of the TASEP is the trivial one.

For our system  introduced above, the Markov matrices $M_L$ are of size $2^L$. 
As mentioned, the entries of these matrices are indexed w.r.t. the
naturally ordered basis of binary vectors of length $L$.   
For convenience, here are the first three of these matrices:
\be \label{mmic1}
M_1 = \begin{bmatrix}
  -\alpha & \alpha+\beta \\
  \alpha & -\alpha-\beta
  \end{bmatrix},~~
M_2 = \begin{bmatrix}
  \star 	& \beta 	& \alpha 	& 1 \\
  0        	& \star 	& 1 		& \alpha \\
  \alpha 	& 0 		& \star 	& \beta \\
  0 		& \alpha 	& 0 		& \star 
  \end{bmatrix},~
M_3 = \begin{bmatrix}
   \star	& \beta	& 0 		& 1 		&\alpha	& 0 		& 1 	& 0 \\
   0 		& \star 	& 1		& 0 		& 0 		& \alpha	& 0 	& 1 \\
   0 		& 0 		& \star	& \beta 	& 1 		& 0 		& \alpha 	& 0 \\
   0 	 	& 0 		& 0 		& \star	& 0 		& 1		& 0 		& \alpha \\
   \alpha 	& 0		& 0		& 0 		& \star	& \beta 	& 0 		& 1 \\
   0 		& \alpha 	& 0 		& 0 		& 0 		& \star	& 1 		& 0 \\
  0 		& 0 		& \alpha  	& 0 		& 0 		& 0 		& \star 	& \beta \\
   0 		& 0 		& 0 		& \alpha	& 0 		& 0 		& 0 		& \star
  \end{bmatrix}.
\ee
As for the diagonal elements $\star$, they have to be set such that the
column sums vanish.

We now state without proof some important results on the Markov matrices 
of the system. These are proved in \cite{am}. We
first show that the Markov matrix itself satisfies a recursion of
order one.

\begin{thm} \label{thm:mm}
Let $\sigma$ denote the matrix
$
 \sigma = \begin{bmatrix}
  0 & 1 \\
  1 & 0
  \end{bmatrix},
$
and $\mathbbm{1}_L$ denote the identity matrix of size $2^L$. Then 
\be \label{mmdecomp}
M_L = \left[ \begin{array}{c | c}
M_{L-1} - \alpha (\sigma \otimes \id{L-2}) & \alpha \id{L-1} + (\sigma
\otimes \id{L-2}) \\ 
\hline
\alpha \id{L-1} & M_{L-1} - \id{L-1} - \alpha (\sigma \otimes \id{L-2} )
\end{array} \right],
\ee
where $M_L$ is written as a $2 \times 2$ block matrix with each block
made up of matrices of size $2^{L-1}$.  
\end{thm}

The transfer matrices can also be explicitly constructed by a recursion
of order one.

\begin{thm} \label{thm:tm}
There exist transfer matrices for the model. 
If one writes the transfer matrix from size $2^{L-1}$ to size $2^L$ by a block
decomposition of matrices of size $2^{L-1} \times 2^{L-1}$ as
\be \label{tmdecomp}
T_{L-1,L} = \left[ \begin{array}{c}
T_1^{(L-1)} \\
\hline
T_2^{(L-1)}
\end{array} \right], 
~~\text{then the matrix $T_{L,L+1}$ can be written as}~~
T_{L,L+1} =   \left[ \begin{array}{c}
T_1^{(L)} \\
\hline
T_2^{(L)}
\end{array} \right],  with 
\ee
\be \label{conjtm1}
 T_1^{(L)} =   
 	\left[ \begin{array}{c|c}
	\ds T_{1}^{(L-1)} + {\alpha}^{-1} T_2^{(L-1)} \,\,  & \,\, 
 	 \ds 2 T_2^{(L-1)} +  {\alpha}^{-1}  T_2^{(L-1)} \\
\\[-3mm]
\hline
\\[-3mm]
	(\sigma\otimes \id{L-2})
	\ds T_2^{(L-1)}\,\,    & \,\,   \ds  {\alpha}^{-1}  T_2^{(L-1)}
	\end{array} \right],  \,\, 
T_2^{(L)} =  \left[ \begin{array}{c|c}
\ds 2T_2^{(L-1)}\,\,   & \,\,    \ds T_2^{(L-1)}  (\sigma\otimes
\id{L-2}) \\ 
\\[-3mm]
\hline
\\[-3mm]
\ds 0 \,\,     &\,\,   T_2^{(L-1)}
\end{array} \right].
\ee
This, along with the initial condition
\be \label{tmic}
T_{1,2} = \left[ \begin{array}{c c}
1 + \beta + \alpha \beta & \alpha + \beta + \alpha \beta \\
\alpha & 1 \\
\alpha + \alpha \beta & \alpha \beta \\
0 & \alpha
\end{array} \right],
\ee
determines recursively  a family of transfer matrices for the matrices $M_L$.
\end{thm}

We can also use the transfer matrices to calculate properties of the
steady state distribution of the Markov process. One quantity of interest is the
so called normalization factor or partition function.

\begin{defn} \label{def:z}
Let the entries of the kernel $\ket{v_L}$ of $M_L$ be normalized 
so that their sum is 1 and each entry written in rationally reduced form.
Then the  {\em partition function}  $Z_L$ for the system of size $L$ is the least
common multiple of the denominators of the entries of  $\ket{v_L}$.
\end{defn}

Because of the way the transfer matrix has been constructed $Z_L$  is
the sum of the 
  entries in  $v_L$. 
For example, the system of size one has
$
\ket{v_1} = \begin{bmatrix}
\alpha + \beta\\
\alpha 
\end{bmatrix},
$
whence $Z_1 = 2\alpha + \beta$. 

\begin{cor} \label{cor:den}
The partition function of the system of size $L$ is
given by
\be  \label{zeq}
Z_L = 2^{\binom{L-1}{2}} (1+2 \alpha)^{L-1} (1+\beta)^{L-1} (2 \alpha
+ \beta).
\ee
\end{cor}

\section{Spectrum of the Markov matrices}
 \label{sec:mm}

 In this section,  we consider the eigenvalues of the Markov matrices
 $M_L$ of the asymmetric annihilation process. 
The following result was stated as a conjecture in
\cite{am}. This will be a corollary of the main result (Theorem \ref{thm:main}) of this
article. 

\begin{thm}\label{thm:conj}
Let the polynomials $A_L(x)$ and $B_L(x)$ be defined as
\be
A_L(x) = \prod_{k=0}^{\lceil L/2 \rceil} (x+2k)^{\binom{L-1}{2k}}, \hskip5mm
B_L(x) = \prod_{k=0}^{\lfloor L/2 \rfloor}
(x+2k+1)^{\binom{L-1}{2k+1}}.
\ee
Then the characteristic polynomial $P_L(x)$ of $M_L$ is given by 
\be \label{charpoly}
P_L(x) = A_L(x) A_L(x+2\alpha+\beta) B_L(x+\beta) B_L(x+2\alpha),
\ee
and successive ratios of characteristic polynomials are given by
\be
\frac{P_{L+1}(x)}{P_L(x)} = B_L(x+1) B_L(x+2\alpha+\beta+1) A_L(x+\beta+1)
A_L(x+2\alpha+1). 
\ee
\end{thm}
 
 This gives only $2L$ distinct eigenvalues out of a possible
$2^L$. There is therefore the question of diagonalizability of the
Markov matrix. Ayyer and Mallick \cite{am} further conjecture the
following.
 
\begin{conj}\label{theconj}
The matrix $M_L$ is maximally degenerate in the sense
that it has exactly $2L$ eigenvectors.
\end{conj}

For $L \geq 1$ we regard $\mathbb{B}^L$ as the vector space  
 	of bitvectors of length $L$ (over the binary field) . The usual scalar product
	of vectors $\bs{b},\bs{c} \in \mathbb{B}^L$ will be denoted by $\bs{b}\cdot\bs{c}$.		
	We will take the set 
	$
	V_L = \left\{ \ket{\bs{b}}\,;\,\bs{b} \in \mathbb{B}^L \right\}
	$
	as  the standard basis of a $2^L$-dimensional (real or complex) vector space, which
	we denote by $\mathcal{V}_L$. Indeed, we will consider $\mathcal{V}_L$ as a vector space
	over an extension over the real or complex field which contains all the
	variables that we introduce below. To be precise, we take
	$\mathcal{V}_L$ as a vector space over a field of rational functions which extends
	the real or complex field. 
	
	The following definitions of linear transformations,
	when considered as matrices, refer to this basis, if not stated otherwise.
	$\mathcal{V}_L$ is the $L$-th tensor power of the 2-dimensional space 
	$\mathcal{V}_1$ in an obvious way.

The transformation $\sigma$ of $\mathcal{V}_1$ is given by the matrix
	$
	\sigma = \begin{bmatrix} 0 & 1 \cr 1 & 0 \end{bmatrix}
	$
	and this extends naturally to transformations $\sigma^{\bs{b}}$ of $\mathcal{V}_L$
	for $\bs{b}= b_1 b_2 \ldots b_L \in \mathbb{B}^L$:
	\[
	\sigma^{\bs{b}}=
	\sigma^{b_1 b_2 \ldots b_L} = 
	\sigma^{b_1} \otimes \sigma^{b_2}   \otimes  \cdots  \otimes \sigma^{b_L}.
	\]

\begin{defn}
	For a vector $\bs{\alpha}=(\alpha_{\bs{b}})_{\bs{b} \in \mathbb{B}^L}$  of variables we define
	the transformation $\mathcal{A}_L(\bs{\alpha})$ of $\mathcal{V}_L$ as
	\[
	A_L(\bs{\alpha}) = \sum_{\bs{b} \in \mathbb{B}^L} \alpha_{\bs{b}}\,\sigma^{\bs{b}}.
	\]
\end{defn}

 A direct way to define these matrices is
	$
	\langle \, \bs{b} \,| \, A_L \,|\, \bs{c} \, \rangle = 
	  \alpha_{\bs{b} \oplus \bs{c}},~(\bs{b},\bs{c} \in \mathbb{B}^L),
	$
	where $\oplus$ denotes the component wise mod-2-addition (exor) of bit vectors.

 For $1 \leq j \leq L$ we define the involutive mappings
	\[
	\phi_j : \mathbb{B}^L \rightarrow \mathbb{B}^L :
	b_1 \ldots   \,b_{j-1}b_{j} \,  b_{j+1}\ldots b_L \mapsto \phi_j \bs{b}= 
	 b_1 \ldots b_{j-1}\, \overline{b_{j}}\,  b_{j+1}   \ldots b_L
	 \]
	 by complementing the $j$-th component, and involutions
	 \[
	 \psi_j : \mathbb{B}^L \rightarrow \mathbb{B}^L : \bs{b} \mapsto \phi_{j} \phi_{j+1} \bs{b}
	 \]
	 by complementing components indexed $j$ and $j+1$,
	 where $\psi_L$ is the same as $\phi_L$.

\begin{defn}
	\begin{enumerate}
	\item 
	For $1 \leq j \leq L$ we define the projection operators $\mathcal{P}_{L,j}$
	acting on $\mathcal{V}_L$ by
	\[
	\mathcal{P}_{L,j} = \sum_{\bs{b} \in \mathbb{B}^L }
		\ket{\bs{b}}\bra{\bs{b}} - \ket{\psi_j^{b_j}(\bs{b})}\bra{\bs{b}}
	\]
	\item 
	For a vector $\bs{b}=(\beta_1,\beta_2,\ldots,\beta_L)$ of variables we put
	$
	\mathcal{B}_L(\bs{\beta}) = \sum_{1 \leq j \leq L} \beta_j\,\mathcal{P}_{L,j}
	$.\\
	$B_L(\bs{\beta})$ denotes the matrix representing $\mathcal{B}_L(\bs{\beta})$
	in the standard basis $V_L$.
	\end{enumerate}
\end{defn}

	Note that in the sum for $\mathcal{P}_{L,j}$ only summands for which $b_j = 1$, 
	i.e., for which $\psi_j(\bs{b}) < \bs{b}$, occur. 
	Indeed: this condition allows only for two situations to
contribute:
	\begin{align} \label{bj1}
	b_{j}\,b_{j+1} = 10 & ~\mapsto~ \overline{b_{j}}\,\overline{b_{j+1}}= 01 & \text{(right shift)} \\\label{bj2}
	b_{j}\, b_{j+1}= 11 & ~\mapsto~ \overline{b_{j}}\,\overline{b_{j+1}}  = 00 & \text{(annihilation)}
	\end{align} 
	Thus these operators $\psi_j$ encode the transitions of our model.
	Also note that by its very definition $B_L(\bs{\beta})$ is an upper triangular matrix.

Writing $\bs{\beta}=( \beta,\gamma,\delta)$ instead of $(\beta_1,\beta_2,\beta_3)$ we have
	for $L=3$
	\[
	B_3(\beta,\gamma,\delta) =
	\begin{bmatrix}
	0 & -   \delta & 0 & - \gamma & 0 & 0 & - \beta & 0 \cr
	0 &   \delta & -\gamma & 0 & 0 & 0 & 0 & - \beta \cr
	0 & 0 & \gamma & -   \delta & - \beta & 0 & 0 & 0 \cr
	0 & 0 & 0 &   \delta+\gamma& 0 & - \beta & 0 & 0 \cr
	0 & 0 & 0 & 0 &  \beta & -  \delta & 0 & -\gamma \cr
	0 & 0 & 0 & 0 & 0 &   \delta+ \beta & -\gamma & 0 \cr
	0 & 0 & 0 & 0 & 0 & 0 & \gamma+ \beta &  -  \delta \cr
	0 & 0 & 0 & 0 & 0 & 0 & 0 &   \delta+\gamma+ \beta
	\end{bmatrix}.
	\]

Our main concern is now with the transformation given by
	\[
	\mathcal{M}_L(\bs{\alpha},\bs{\beta})= \mathcal{A}_L(\bs{\alpha})-\mathcal{B}_L(\bs{\beta}).
	\]	
	Before we can state the main result we have to introduce some more
	notation, But before doing so, we note that the corresponding matrix
	${M}_L(\bs{\alpha},\bs{\beta})= {A}_L(\bs{\alpha})-{B}_L(\bs{\beta})$ reduces to
	the matrix $M_L$ above when properly specialized:

\begin{lem}\label{firstlemma}
	We have $A_L(\bs{\alpha}')-B_L(\bs{\beta}')=M_L$ for
	 $\bs{\alpha}'=(\alpha_{\bs{b}}')_{\bs{b} \in \mathbb{B}^L}$ 
	and $\bs{\beta}' = ( \beta_j')_{1 \leq j \leq L}$  given by
	\[
		\alpha_{\bs{b}}' =
		\begin{cases} 
		- \alpha & \text{if}~\bs{b} = 00\ldots 00 \cr
		 \alpha&\text{if}~\bs{b} = 10\ldots 00 \cr      
	          0  &\text{otherwise}
	          \end{cases}
	          ~~~\text{and}~~
	          \beta' =
	          \begin{cases}
	          1 & \text{if}~1 \leq j < L \cr
	          \beta &\text{if}~j=L
	          \end{cases}.
	\]
	
\end{lem}

	We will now consider the  transformation
	$\mathcal{M}_L(\bs{\alpha},\bs{\beta})$ in a different basis of $\mathcal{V}_L$.	
	Let
	$
	H = \frac{1}{\sqrt{2}}\begin{bmatrix} 1 & 1 \cr 1 & -1 \end{bmatrix}
	$
	be the familiar Hadamard matrix and define $H_L$ as its
	$L$-th tensor power, the matrix$L$ Hadamard transform  of order $L$:
	\[
	H_L = H^{\otimes L}=\frac{1}{2^{L/2}}\,
	\left[\, (-1)^{\bs{b}\cdot\bs{c}} \, \right]_{\bs{b},\bs{c} \in \mathbb{B}^L}.
	\]
	The columns of this matrix, denoted by
	$\ket{ w^{\bs{b}}} = H_L \ket{\bs{b}}$ for $\bs{b} \in \mathbb{B}^L$,
	form an orthonormal basis 
	$W_L= \left\{ H_L \ket{\bs{b}}\,;\,\bs{b}\in \mathbb{B}^L \right\}$ of $\mathcal{V}_L$.
	 The following assertion is easily checked:
	 
\begin{lem}\label{secondlemma}
	 The (pairwise commuting) transformations $\sigma^{\bs{c}}~(\bs{c} \in \mathbb{B}^L)$ 
	diagonalize in the $W_L$-basis. More precisely:
	\[
	\sigma^{\bs{c}} \, \ket{ w^{\bs{b}}}=(-1)^{\bs{b}\cdot\bs{c}}\,\ket{w^{\bs{b}}}~~~
	(\bs{b},\bs{c} \in \mathbb{B}^L).
	\]
	Thus also the transformation $\mathcal{A}_L$ diagonalizes 
	in the $W_L$-basis and its eigenvalues are given by 
	\[
	(H_L \cdot A_L \cdot H_L)\ket {w^{\bs{b}}}=\lambda_{\bs{b}}\ket {w^{\bs{b}}}
	\]
	where
	$
	\lambda_{\bs{b}} = \sum_{c \in \mathbb{B}^L}\alpha_{\bs{c}}\,(-1)^{\bs{b}\cdot \bs{c}}
	=  \sum_{c \in \mathbb{B}^L}\alpha_{\bs{c}}\,\langle\, \bs{b}\,|\,H \,|\,\bs{c}\,\rangle.
	$
\end{lem}

 The crucial observation is now the following: even though the transformation  
	$\mathcal{A}_L$ diagonalizes in the $W_L$-basis, the transformation $\mathcal{B}_L$ doesn't, 
	it is not even triangular in this basis. 
	But it turns out that a slight modification of the $W_L$-basis will be suitable
	for at the same time diagonalizing  $\mathcal{A}_L$  and bringing the $\mathcal{B}_L$
	in (lower) triangular form. For that purpose we introduce the invertible linear transformation
	\[
	\Delta : \mathbb{B}^L \rightarrow \mathbb{B}^L :
	\bs{b} = b_1b_2 \ldots b_L \mapsto  \bs{b}^\Delta 
	= \left[\sum_{1 \leq i \leq L-j+1}b_i \right]_{1 \leq j \leq L}
	\] 
	where the sum has to be taken in the binary field.
	As an example ($L=3$):
	\[
	\begin{array}{r|cccccccc}
	\bs{b} 		& 000 & 001 &  010 & 011 & 100 & 101& 110 & 111 \cr\hline
	\bs{b}^\Delta 	& 000 & 100 & 110 & 010 & 111 & 011 & 001 & 101
	\end{array}
	\]
 The basis $\widetilde{W}_L = \{ \ket{w^{ \bs{b}^\Delta}} \}$ is nothing but a
	rearrangement of the  $W_L$-basis, hence the transformation $\mathcal{A}_L$ 
	diagonalizes in this basis as well (with the corresponding eigenvalues). We will write
	$\widetilde{H}_L$ for the rearrangement of the Hadamard matrix
	in this new ordering of the elements of $\mathcal{B}^L$.
	The clue is now contained in the following proposition: 

\begin{prop} \label{mainprop}
	\[
	\widetilde{H}_L \cdot B_L(\bs{\beta}) \cdot \widetilde{H}_L = 
	B_L^{\mathtt{t}}(\bs{\beta}^{\mathtt{rev}})
	\]
	where $\bs{\beta}^{\mathtt{rev}}=(\beta_L,\beta_{L-1}, \ldots, \beta_1)$ is the
	reverse of $\bs{\beta} = (\beta_1,\beta_2, \ldots,\beta_L)$,
	 $\mathtt{t}$ denoting  transposition.
\end{prop}

We illustrate this proposition in the case $L=3$ by displaying matrices
$2^{3/2} \widetilde{H}_3$ (left) and 
$\widetilde{H}_3 \cdot B_3(\beta,\gamma,\delta) \cdot \widetilde{H}_3 =
B_3(\delta,\gamma,\beta)^{\mathtt{t}}$ (right). Note that $\widetilde{H}_L$
is symmetric because $\Delta$ (as a matrix) is symmetric.
\small	\[ 
	 \left[ \begin {array}{rrrrrrrr} 1&1&1&1&1&1&1&1\\\noalign{\medskip}1&
	1&1&1&-1&-1&-1&-1\\\noalign{\medskip}1&1&-1&-1&-1&-1&1&1
	\\\noalign{\medskip}1&1&-1&-1&1&1&-1&-1\\\noalign{\medskip}1&-1&-1&1&-
	1&1&1&-1\\\noalign{\medskip}1&-1&-1&1&1&-1&-1&1\\\noalign{\medskip}1&-
	1&1&-1&1&-1&1&-1\\\noalign{\medskip}1&-1&1&-1&-1&1&-1&1\end {array}
	 \right] ,~~~~
	 \left[ \begin {array}{cccccccc} 0&0&0&0&0&0&0&0\\\noalign{\medskip}-
 	\beta& \beta&0&0&0&0&0&0\\\noalign{\medskip}0&-\gamma&\gamma&0&0&0&0&0
	\\\noalign{\medskip}-\gamma&0&- \beta& \beta+\gamma&0&0&0&0
	\\\noalign{\medskip}0&0&-  \delta&0&  \delta&0&0&0\\\noalign{\medskip}0&0&0&
	-  \delta&- \beta& \beta+  \delta&0&0\\\noalign{\medskip}-  \delta&0&0&0&0&-
	\gamma&\gamma+  \delta&0\\\noalign{\medskip}0&-  \delta&0&0&-\gamma&0&-
	 \beta& \beta+  \gamma+\delta\end {array} \right] .
	\]
\normalsize
	 
 The proof of Proposition \ref{mainprop} will be given below. It leads 
	to the main result by looking at the matrix representation
	of $\mathcal{M}_L(\bs{\alpha},\bs{\beta})=\mathcal{A}_L-\mathcal{B}_L$ in the $\widetilde{W}_L$-basis,
	where it takes lower triangular form. Hence the eigenvalues, which are
	$2^L$ pairwise distinct linear polynomials in the $\alpha$- and $\beta$-variables,
	can be read directly from the main diagonal. In contrast to Theorem  \ref{thm:conj},
	the ex-conjecture, all eigenvalues are simple.
	
\begin{thm}\label{thm:main}
	\[
	\det M_L(\bs{\alpha},\bs{\beta})= 
	\det \left[ A_L(\bs{\alpha}) - B_L(\bs{\beta})\right] = 
	\prod_{\bs{b} \in \mathbb{B}^L} \left(
	\lambda_{\bs{b}^\Delta} -
	  \bs{\beta}^{\mathtt{rev}} \cdot  \bs{b}  \right)
	\]
\end{thm}

Illustration of the Theorem for $L=3$ recalling $\bs{\beta}=(
\beta,\gamma,\delta) =(\beta_1,\beta_2,\beta_3)$:
	\[
	\begin{array}{cccl}
	\bs{b}  & \bs{b}^\Delta  & \lambda_{b^\Delta} & (\delta,\gamma,\beta)\cdot\bs{b} \cr
	000 & 000 & [++++++++]\cdot \bs{\alpha} & 0 \cr
	001 & 100 & [++++----]\cdot \bs{\alpha}     & \beta \cr
	010 & 110 & [++----++]\cdot \bs{\alpha}     & \gamma \cr
	011 & 010 & [++--++--]\cdot \bs{\alpha}    & \beta+\gamma  \cr
	100 & 111 & [+--+-++-]\cdot \bs{\alpha}     &\delta \cr
	101 & 011 & [+--++--+]\cdot \bs{\alpha}     &\beta+\delta \cr
	110 & 001 & [+-+-+-+-]\cdot \bs{\alpha}     &  \gamma+\delta \cr
	111 & 101 & [+-+--+-+]\cdot \bs{\alpha}     & \beta+\gamma+\delta
	\end{array}
	\]

So, as an example, the line for $\bs{b}=101$ contributes the factor
$\alpha_{000}-\alpha_{001}-\alpha_{010}-\alpha_{011}
+\alpha_{100}-\alpha_{101}-\alpha_{110}-\alpha_{111} - \beta-\delta$ to the product.

	To prepare for the proof of  Proposition \ref{mainprop}
	we state without proof
	simple relations between the transformations $\psi_j, \phi_{L-j+1}$ and $\Delta$:
	 
\begin{lem}\label{thirdlemma} For $\bs{b},\bs{c} \in \mathbb{B}^L$ and $1 \leq j \leq L$ we have
	\begin{enumerate}
	\item
	$
	( \psi_j \bs{b})^\Delta  = \phi_{L-j+1}(\bs{b}^{\Delta} )
	$
	\item
	$
	\bs{b}^{\Delta} \cdot \psi_j \bs{c} = \bs{b}^\Delta \cdot \bs{c} + b_{L-j+1}
	$
	 \end{enumerate}
\end{lem}	 
	 Fact 2. is a consequence of fact 1.
	 
{\medskip\noindent\bf Proof of Proposition \ref{mainprop}}.~
	The actions of the transformations $\mathcal{P}_{L,j}$, seen
	in the $\widetilde{W}_L$-basis, are given by:
	\[
	\mathcal{P}_{L,j} : \ket{w^{\bs{b}^\Delta}} \mapsto 
	\begin{cases}
	- \ket{w^{(\psi_{L-j+1} \bs{b})^{\Delta}}} & \text{if}~~b_{L-j+1}=0, \cr
	\ket{w^{\bs{b}^\Delta}}&  \text{if}~~ b_{L-j+1}=1.
	\end{cases}
	\]
	To see this, we compute
	\begin{align*}
	\mathcal{P}_{L,j}\ket{w^{\bs{b}^\Delta}} 
	&=\sum_{\bs{c}\in\mathbb{B}^L}
	\langle\,\bs{b}^\Delta\,|\,H\,|\,\bs{c}\,\rangle\mathcal{P}_{L,j} |\,\bs{c}\,\rangle
	=\sum_{\bs{c}>\psi_j\bs{c}}
	\langle\,\bs{b}^\Delta\,|\,H\,|\,\bs{c}\,\rangle
	\left(|\,\bs{c}\,\rangle-|\,\psi_j\bs{c}\,\rangle\right)\cr
	&=\sum_{\bs{c}>\psi_j\bs{c}}
	\langle\,\bs{b}^\Delta\,|\,H\,|\ \bs{c}\,\rangle\,|\,\bs{c}\,\rangle-\sum_{\bs{c}<\psi_j\bs{c}}
	\langle\,\bs{b}^\Delta\,|\,H\,|\,\psi_j\bs{c}\,\rangle|\,\bs{c}\,\rangle,
	\end{align*}
	using the involutive nature of $\psi_j$ for the second sum. 
	Now, using 2. from Lemma \ref{thirdlemma},
	\begin{align*}
	\langle\, \bs{b}^\Delta \,|\, H \, | \, \psi_j \bs{c} \, \rangle =
	(-1)^{\bs{b}^\Delta\cdot \psi_j \bs{c}} = (-1)^{\bs{b}^\Delta\cdot\bs{c}+b_{L-j+1}}
	= (-1)^{b_{L-j+1}} \langle\, \bs{b}^\Delta \,|\, H \, | \,  \bs{c} \, \rangle
	\end{align*}
	and thus
	\[
	\mathcal{P}_{L,j} \ket{ w^{\bs{b}^\Delta}} =
	\sum_{\bs{c} > \psi_j \bs{c}}
		\langle\, \bs{b}^\Delta \,|\, H \, | \, \bs{c} \, \rangle \,|\,\bs{c} \, \rangle - (-1)^{b_{L-j+1}}
		\sum_{\bs{c} < \psi_j \bs{c}}
		\langle\, \bs{b}^\Delta \,|\, H \, | \,  \bs{c} \, \rangle |\,\bs{c} \, \rangle.
	\]
	The conclusion in the case $b_{L-j+1}=1$ is now obvious.
	
	As for the  case $b_{L-j+1}=1$, we see, using item 1. from Lemma \ref{thirdlemma}, that
	\begin{align*}
	\ket{w^{(\psi_{L-j+1} \bs{b})^{\Delta}}} &=
		\sum_{\bs{c} \in \mathbb{B}^L} 
			\langle\, (\psi_{L-j+1} \bs{b})^\Delta \,|\, H \, | \, \bs{c} \, \rangle \, |\,\bs{c} \, \rangle
		=
		\sum_{\bs{c} \in \mathbb{B}^L} 
			(-1)^{(\psi_{L-j+1} \bs{b})^\Delta \cdot \bs{c}} \, |\,\bs{c} \, \rangle \cr
		&=
		\sum_{\bs{c} \in \mathbb{B}^L} 
			(-1)^{\phi_{j}(\bs{b} ^\Delta) \cdot \bs{c}} \, |\,\bs{c} \, \rangle 
		=
		\sum_{\bs{c} \in \mathbb{B}^L} 
			(-1)^{ \bs{b}^\Delta \cdot \bs{c} + c_{j}}\,  |\,\bs{c} \, \rangle \cr
		&=
		- \sum_{\bs{c}: c_j=1} \langle\, \bs{b}^\Delta \,|\, H \, | \, \bs{c} \, \rangle \,|\,\bs{c} \, \rangle
		+ \sum_{\bs{c}: c_j=0} \langle\, \bs{b}^\Delta \,|\, H \, | \, \bs{c} \, \rangle \,|\,\bs{c} \, \rangle\
		= - \mathcal{P}_{L,j} \ket{ w^{\bs{b}^\Delta}}.\hfill~~~~~~~ \Box
	\end{align*}

\begin{cor} If we consider the special case where
	$\alpha_{\bs{b}}=0$ for all $\bs{b} \in \mathbb{B}^L$, except 
	$\alpha_{00\ldots 0}= \alpha_0$ and $\alpha_{10\ldots00}=\alpha_1$, 
	and where $\beta_1=\ldots=\beta_{L-1}=1$ and $\beta_L=\beta$, then
	the determinant of the Theorem simplifies to the product 
	$\Pi_1 \cdot \Pi_2 \cdot \Pi_3 \cdot \Pi_4$ of the following four terms:
	\begin{align*}
	\Pi_1 &=\prod_{0 \leq 2k < L} (\alpha_0+\alpha_1-2k)^{\binom{L-1}{2k}} &
	\Pi_2 &=\prod_{0 \leq 2k-1 < L} (\alpha_0+\alpha_1- \beta -2k+1)^{\binom{L-1}{2k-1}} \cr
	\Pi_3 &=\prod_{0 \leq 2k-1 < L} (\alpha_0-\alpha_1  -2k+1)^{\binom{L-1}{2k-1}} &
	\Pi_4 &=\prod_{0 \leq 2k < L} (\alpha_0-\alpha_1-\beta-2k)^{\binom{L-1}{2k}} 
	\end{align*}
\end{cor}
	
	For the proof note that  each $\bs{b} \in \mathbb{B}^L$ we get  as the contribution from
	$\mathcal{A}_L$
	\begin{align*}
	\lambda_{\bs{b}^{\Delta}}\cdot\bs{\alpha} &= 
\sum_{\bs{c} \in \mathbb{B}^L}\alpha_{\bs{c}}\langle\,\bs{b}^{\Delta}\,|\,H\,|\,\bs{c}\,\rangle\cr
		&=
		\alpha_0\,\langle\bs{b}^{\Delta}\,|\,H\,|\,{00\ldots00}\rangle+
		\alpha_1\,\langle\bs{b}^{\Delta}\,|\,H\,|\,{10\ldots00}\rangle\cr
		&= 
		\alpha_0(-1)^{\bs{b}\cdot\Delta\cdot 00\ldots 00}+
		\alpha_1(-1)^{\bs{b}\cdot\Delta\cdot 10\ldots 00}
		=
		\alpha_0+(-1)^{\|\bs{b}\|}\,\alpha_1
		\end{align*}
	because $\Delta \cdot 00 \ldots 00=11 \ldots 11$ and then 
	$\bs{b}\cdot 11\ldots 11 \equiv \|\bs{b}\| \bmod 2$, where  $\|\bs{b}\|$
	denotes the Hamming weight of $\bs{b}$ and where we have used the fact that $\Delta$
	is a symmetric matrix. Thus the $2^L$ eigenvalues are
	\[
	\alpha_0  + (-1)^{\| \bs{b} \|}\, \alpha_1 - \bs{\beta}^{\mathtt{rev}} \cdot \bs{b}~~~(\bs{b} \in \mathbb{B}^L).
	\]
	Now there are four cases to consider:
	\begin{enumerate}
	\item $\| \bs{b} \|$ is even and $b_1=0$: this gives eigenvalues
		$
		\alpha_0+\alpha_1 - \| b_2 b_3 \ldots b_L\|
		$
		and since $b_1$ does not contribute to $\|\bs{b}\|$ the 
		vector $b_2 b_3 \ldots b_L$ must have even weight $2k$.
		There are $\binom{L-1}{2k}$ possibilities which account for $\Pi_1$.
	\item $\| \bs{b} \|$ is even and $b_1=1$: this gives eigenvalues
		$
		\alpha_0+\alpha_1 - \beta - \| b_2 b_3 \ldots b_L\|
		$
		and since $b_1$ does  contribute to $\|\bs{b}\|$ the 
		vector $b_2 b_3 \ldots b_L$ must have odd weight $2k-1$.
		There are $\binom{L-1}{2k-1}$ possibilities  which account for $\Pi_2$.
	\item $\| \bs{b} \|$ is odd and $b_1=0$: this gives eigenvalues
		$
		\alpha_0-\alpha_1 - \| b_2 b_3 \ldots b_L\|
		$
		and since $b_1$ does not contribute to $\|\bs{b}\|$ the 
		vector $b_2 b_3 \ldots b_L$ must have  odd weight $2k-1$.
		There are $\binom{L-1}{2k-1}$ possibilities  which account for $\Pi_3$.
	\item $\| \bs{b} \|$ is odd and $b_1=1$: this gives eigenvalues
		$
		\alpha_0-\alpha_1 - \beta - \| b_2 b_3 \ldots b_L\|
		$
		and since $b_1$  does contribute to $\|\bs{b}\|$ the 
		vector $b_2 b_3 \ldots b_L$ must have even weight $2k$.
		There are $\binom{L-1}{2k}$ possibilities which account for $\Pi_4$.
	\end{enumerate}
	
\begin{cor}\label{cor:spec}
	 Setting now $\alpha_0=-\alpha$ and $\alpha_1=\alpha$, 
	i.e., specializing as in Lemma \ref{firstlemma}, gives for 
	$\det M_L$ 
	a product $\Pi_1' \cdot \Pi_2' \cdot \Pi_3' \cdot \Pi_4'$ of the following four terms:
	\begin{align*}
	\Pi_1' &=\prod_{0 \leq 2k < L} (-2k)^{\binom{L-1}{2k}} &
	\Pi_2 &=\prod_{0 \leq 2k-1 < L} (- \beta -2k+1)^{\binom{L-1}{2k-1}} \cr
	\Pi_3 &=\prod_{0 \leq 2k-1 < L} (-2 \alpha  -2k+1)^{\binom{L-1}{2k-1}} &
	\Pi_4 &=\prod_{0 \leq 2k < L} ( -2\alpha_1-\beta-2k)^{\binom{L-1}{2k}} 
	\end{align*}
	which are precisely the $2\,L$ distinct eigenvalues of the original Conjecture.
\end{cor}

\section{Concluding remarks}

We have been able to solve Ayyer and Mallick's conjecture about the
eigenvalues of the asymmetric annihilation process by embedding it
into a more general model and using an orthogonal transform which
makes the transition matrix upper triangular.  In sharp contrast to the
original problem, the general situation with parameters
$\alpha_{\bs{b}}$ and $b_j$ 
(which may be given a ``physical'' interpretation using
$\langle \, \bs{b} \,| \, A_L \,|\, \bs{c} \, \rangle = 
	  \alpha_{\bs{b} \oplus \bs{c}}$ and \eqref{bj1},\eqref{bj2}) 
is easier to handle because it is not
degenerate: all ``symbolic'' eigenvalues are simple. Our proof does
not seem to explain the maximum amount of degeneracy, as stated in
Conjecture \ref{theconj}.
	
On the other hand, we mention that the result of Corollary
\ref{cor:den} about the partition function can be extended to the more
general model.  Again, in contrast to the inductive approach of Ayyer
and Mallick in \cite{am}, as outlined in Sec. \ref{sec:tma}, we can
solve this problem directly by transforming it orthogonally into the
basis where it shows its triangular structure.
	
We start by remarking that the columns sums of the extended model are
constant $\sum \alpha_c$, though not zero. This implies that
$\bra{1_L}$ is the unique left eigenvector with eigenvalue
$\overline{\alpha}= \sum \alpha_c$ of $M_L(\alpha,\beta)$.  The right
eigenvector $\ket{\bs{x}}$ with the same eigenvalue corresponds to the
steady state distribution of the original problem. Then $\ket{\bs{y}}
= \widetilde{H}_L \ket{\bs{x}}$ satisfies
$\widetilde{M}_L(\bs{\alpha},\bs{\beta}) \ket{\bs{y}}=
\overline{\alpha}\ket{\bs{y}}$, where $\widetilde{M}_L =
\widetilde{H}_L \cdot M_L(\bs{\alpha},\bs{\beta}) \cdot
\widetilde{H}_L$ is the matrix seen in the
$\widetilde{W}_L$-basis. This triangular system for $\bs{y}$ is
written explicitly as \[ \left( \lambda_{b^\Delta}^* +
\bs{\beta}^{\mathtt{rev}} \cdot \bs{b} \right)\, y_{\bs{b}} =
\sum_{j\,:\,b_j=1} \beta_{L-j+1} y_{\psi_j \bs{b}}~~~~(\bs{b} \in
\mathbb{B}^L) \] where \[ \lambda_{\bs{b}}^* = \overline{\alpha}-
\lambda_{\bs{b}} = 2 \sum_{\bs{c}\,:\,\bs{b}\cdot
\bs{c}=1}\alpha_{\bs{c}}~~~~~(\bs{b} \in \mathbb{B}^L).  \] Note that
the sum on the right only contains terms $y_{\bs{c}}$ where
$\bs{c}=\psi_j \bs{b} < \bs{b}$.  For $\bs{b}=00\ldots0$ the equation
is void, so we may put $y_{00\ldots 0}=1$.  Since the polynomials
$\lambda_{b^\Delta}^* + \bs{\beta}^{\mathtt{rev}} \cdot \bs{b}$ are
mutually coprime, this shows by induction that the denominator of the
rational normal form of $y_{\bs{b}}$ is the product of all polynomials
$\lambda_{c^\Delta}^* + \bs{\beta}^{\mathtt{rev}} \cdot \bs{c}$, where
$\bs{c}$ runs over the binary vectors that can be obtained from
$\bs{b}$ by successive application of decreasing
$\psi_j$-transformations.  Consequently, the product of linear
polynomials \[ Z(\bs{\alpha},\bs{\beta})=
\prod_{\bs{0}\neq\bs{b}\in\mathbb{B}^L} \left(\lambda_{b^\Delta}^* +
\bs{\beta}^{\mathtt{rev}} \cdot \bs{b}\right) \] is the least common
multiple of the denominators of the $y_{\bs{b}}$.  This property is
invariant under the Hadamard transform, so it applies also to the
coefficients of $\ket{\bs{x}} = \widetilde{H}_L \ket{\bs{y}}$. But \[
\langle \bs{1}_L | \bs{x} \rangle = \langle \bs{1}_L | \widetilde{H}_L
| \bs{y} \rangle = 2^{L/2} \, \langle 100 \ldots 00 | \bs{y} \rangle =
2^{L/2} \, y_{00\ldots 0}=2^{L/2}, \] so $2^{-L/2}\, \ket{\bs{x}}$ is
already normalized and can be seen as the ``symbolic'' stationary
distribution in the generalized model. What we have shown is:
	
\begin{thm}\label{thm:partition}
$Z(\bs{\alpha},\bs{\beta})$ is the partition function related to
$M_L(\bs{\alpha},\bs{\beta})$.  
\end{thm}
	
We conclude by remarking that the specialization as in Lemma
\ref{firstlemma} and Corollary \ref{cor:spec} brings us back to
Corollary \ref{cor:den}. This is not completely obvious, since the
expression in Corollary \ref{cor:den} has only $\binom{L+1}{2}$
factors, whereas in Theorem \ref{thm:partition} there are $2^L-1$
factors.  What happens is that upon specialization the requirements
for least common multiples and greatest common divisors change. Taking
this into account one finds that from the general expression for
$Z(\bs{\alpha},\bs{\beta})$ only the $\binom{L+1}{2}$ terms where
$\bs{b} \in \mathbb{B}^L$ with $\|\bs{b}\|=1~\text{or}~2$ contribute
-- and this is precisely the statement of Corollary \ref{cor:den}.

\end{document}